\documentclass{article}

\usepackage[utf8]{inputenc}
\usepackage{amsmath,amsthm,amssymb}
\usepackage{color,soul}
\usepackage{graphicx}
\usepackage{mathrsfs}
\usepackage{multicol}
\usepackage{bbm}
\usepackage{bm}
\usepackage{amsthm}
\usepackage{fullpage}
\usepackage[ruled]{algorithm2e}
\usepackage[colorlinks = true, citecolor = blue, urlcolor = blue]{hyperref}

\newtheorem{lem}{Lemma}

\newtheorem{defn}[lem]{Definition}

\newtheorem{example}{Example}
\newtheorem{rmk}{Remark}

\newcommand{\R}{\mathbb{R}}				    	        
\newcommand{\X}{\mathbb{X}}                             
\newcommand{\U}{\mathcal{U}}                            
\newcommand{\nerve}{\mathcal{N}}                        
\newcommand{\goto}{\rightarrow}			    	        
\newcommand{\LP}{\left(} \newcommand{\RP}{\right)}	    
\newcommand{\LC}{\left\{} \newcommand{\RC}{\right\}}	
\newcommand{\Ln}{\left|} \newcommand{\Rn}{\right|}	    

\providecommand{\keywords}[1]
{
  \small	
  \textbf{\textit{Keywords---}} #1
}

\begin{document}

\title{Multiscale Geometric Data Analysis via Laplacian Eigenvector Cascading \thanks{Support: NSF (DMS-1622301) and  DARPA (HR0011-16-2-003)}}

\author{Joshua L. Mike \\
\textit{Department of Computational Mathematics,} \\
\textit{Science, and Engineering}  \\
\textit{Michigan State University}\\
East Lansing, USA \\
\texttt{mikejosh@msu.edu}
\and
Jose A. Perea \\
\textit{Department of Computational Mathematics,} \\
\textit{Science, and Engineering} \\
\textit{Department of Mathematics} \\
\textit{Michigan State University}\\
East Lansing, USA \\
\texttt{joperea@msu.edu}
}

\maketitle


\begin{abstract}
    We develop here an algorithmic framework for constructing consistent multiscale Laplacian eigenfunctions (vectors) on data.
    Consequently, we address the unsupervised machine learning task of finding scalar functions capturing consistent structure across scales in data,
    in a way that encodes intrinsic geometric and topological features.
    This is accomplished by two algorithms for eigenvector cascading.
    We show via examples that cascading accelerates the computation of graph Laplacian eigenvectors, and more importantly,  
    that one obtains consistent bases of the associated eigenspaces across scales.
    Finally, we present an application to  TDA mapper, showing that our  multiscale Laplacian eigenvectors identify stable flair-like structures in mapper graphs of varying granularity.
\end{abstract}


\section{Introduction}


Metric data is often represented as a graph, with data points as nodes and edges drawn according to comparison predicates.
By adding edge similarity weights one can define graph Laplacian operators, whose eigenvectors and eigenvalues have been shown to encode the underlying geometry/topology of the dataset \cite{Lap_Shape_Rep, Geo_DeRham}.
The result is a collection of real-valued functions on the nodes of the graph---the Laplacian eigenvectors---which can
then be used in tasks like dimensionality reduction \cite{isomap, Lap_Eigen_Map}, semi-supervised learning \cite{Belk_Niyogi_SSLearn}, ranking problems \cite{HodgeRank}, and direct quantification of geometry \cite{Lap_Shape_Match, Lap_Shape_Rep, Lap_Fluid_Sim}.
In particular, null eigenvectors correspond to connected components and eigenvectors with small eigenvalue code for weakly connected components.

Along these lines, here \emph{we use  ideas from persistent cohomology and Hodge theory to construct and relate graphs from a dataset viewed at multiple scales,
 and to make consistent choices of eigenvector features across coarseness levels}.
The multiscale graphs are obtained via a series of refinements on an initial cover of the dataset, and
the resulting graph relationships (simplicial maps) are used to define a pair of eigenvector cascade methods: first cascade (Alg. \ref{alg:cascade1}) and second cascade (Alg. \ref{alg:cascade2}).
We show initially that first cascade speeds up Laplacian Eigenvector computations, while further experiments demonstrate that second cascade reliably tracks weakly connected components across scale;
This can be seen as a soft version of 0-dimensional persistent homology, which tracks connected components in a filtered space. 
By choosing geometrically consistent basis vectors across scales, second cascade in particular overcomes eigenvector instability for clustered or repeated eigenvalues \cite{EV_repeats}.

The proposed methodology is instantiated by two sources of towers of covers:
The cover tree construction (Def. \ref{def:CT}), and TDA mapper \cite{mapper}.
We use the cover tree construction on low-dimensional examples to demonstrate the benefits of the cascade algorithms, and then apply them to multiscale TDA mapper graphs as an analysis application.
Within this context, Laplacian eigenvectors highlight weakly connected components of the mapper graphs and thereby track geometrically conspicuous regions of a dataset and assess their stability to scale.
It is worth mentioning that \cite{MSmapper} also considers  a multiscale notion of mapper.
While here we investigate the persistence of Laplacian eigenspace features, \cite{MSmapper} instead explores its persistent homology.

\subsection{Prior Work}
Large datasets lead to large graphs and Laplacian matrices with prohibitive computations.
Thus, several methods view a graph from multiple scales to accelerate Laplacian eigenvector estimation  \cite{multi-scale_spectral_graph}, \cite{F1_cascade}, \cite{ACE_Cascade}.
Analogously to our approach, these use coarsenings of a large graph and propagate eigenvectors on coarse graphs as initial states for solvers on finer ones.
The graph collapse method in the latter two works is similar in spirit to ours, 
but the topological approach presented here yields functions defined on an open neighborhood of the data, and 
it is thus  amenable to sparsification methods like
landmark subsampling.


Moreover,  we are not simply interested in accelerating computations.
Indeed, our primary goal is to capture the progression of Laplacian eigenspaces and to identify stable (i.e, persistent) features.
This allows a multi-scale view of dataset geometry via important structure which persists across various scales.
The interpretation of Laplacian eigenvectors in terms of 
Hodge theory makes clear the connection to (persistent) cohomology: null eigenvectors of the Hodge Laplacian correspond to unique cohomology classes.
In this fashion, one may also directly compare our methodology with the notion of persistent (co)homology \cite{TDA} which tracks the persistence of (co)homological generators
across scales,  or for a self map \cite{edelsbrunner2015persistent} (further details given in Rmk. \ref{rmk:topo_compare}).

For large eigenvector problems, solving directly is intractable and iterative methods such as  Lanczos \cite{lanczos} are used instead.
For a  multiscale implementation, we are interested in methods which take advantage of an initial state such as Rayleigh quotient iteration \cite{RQI} or conjugate gradient descent \cite{Conjugate_Gradient}.
In particular, we consider locally optimally block preconditioned conjugate gradient (LOBPCG) \cite{lobpcg} which is well-suited for large-scale symmetric problems by its use of parallel processing along a matrix's sparse structure.
LOBPCG has been shown to be very efficient in practice, despite its incomplete theory \cite{LOBPCG_theory}.

\section{Multiscale Laplacian Eigenfunctions on  Data: A topological perspective} \label{S:bg}
Henceforth we assume that the reader is familiar with the topological notions of open covers, simplicial (e.g., nerve) complexes, and partitions of unity.
Given a dataset $X$ inside a metric space $(\mathbb{X}, d)$, let  $\U = \LC U_j \RC$
be an open cover of $X$, and select a partition of unity $\LC \phi_j \RC$ dominated by $\U$.
This yields  a \emph{nerve map}  defined  in barycentric coordinates by
\begin{align*}
    \phi:\cup \U &\goto \Ln \nerve(\U) \Rn \\
    x &\mapsto \LP \phi_1(x),...,\phi_N(x) \RP
\end{align*}
 into the geometric realization of the nerve complex.

\begin{figure}[htb!]
    \centering
    \includegraphics[width = .5\textwidth]{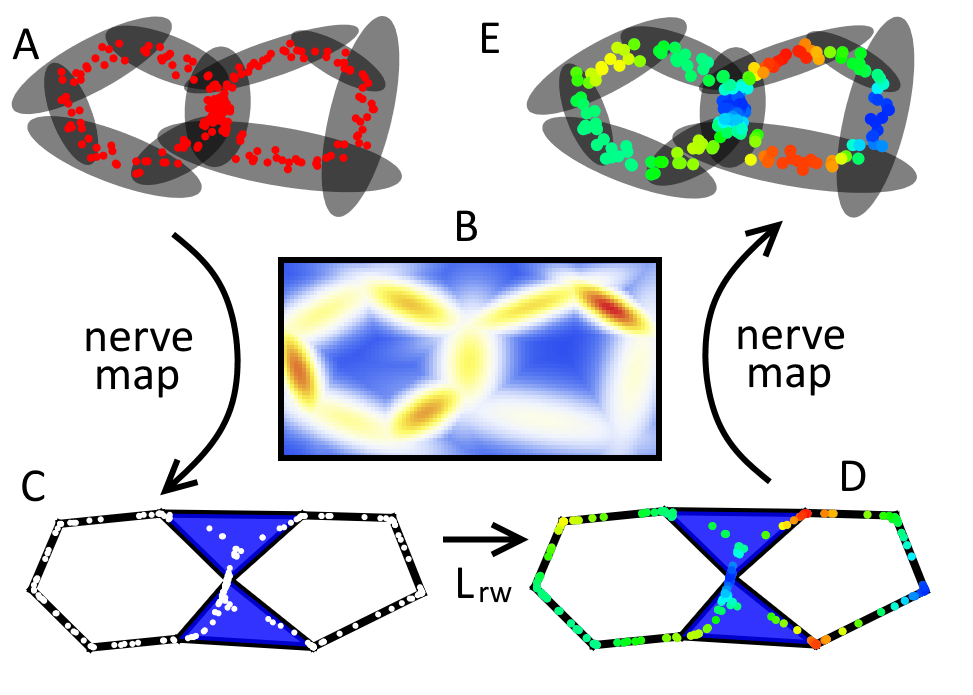}
    \caption{The general pipeline for nerve mapping and defining approximate Laplacian eigenfunctions is shown here.
    (A) a cover $\U$ is provided for the dataset.
    (B) A partition of unity (PoU) subordinate to $\U$ is chosen. 
    (C) The PoU induces a nerve map $\phi:\bigcup\mathcal{U} \goto \Ln\nerve(\U)\Rn$.
    (D) Edge weights induce a Laplacian operator $L_{rw}$, and the associated eigenvectors are scalar functions on the vertices of $\mathcal{N(U)}$, which extend linearly to $|\mathcal{N(U)}|$.
    (E) Precomposing with the nerve map defines approximate Laplacian eigenfunctions on $\cup\mathcal{U}$.}
    \label{fig:NL_vignette}
\end{figure}

The  1-skeleton of $\mathcal{N(U)}$ yields a graph $G$ which describes the data at a fixed scale.
Weighting the edges of this graph yields a Laplacian operator, whose eigenvectors are scalar functions on the vertices of $\nerve(\U)$ and which extend linearly to functions on $|\nerve(\U)|$.
Precomposing the extended eigenvectors with the nerve map $\phi$ yields approximate eigenfunctions on an open neighborhood of the dataset, namely $\cup \U$.
These approximate eigenfunctions can be viewed as geometric coordinates for  $X$, at the scale furnished by the covering $\U$. 
A brief vignette describing the nerve map and associated eigenfunctions is shown in Fig. \ref{fig:NL_vignette}.

By repeatedly refining a cover $\U$ of  $X$, we obtain a tower of covers $\LC \U^i \RC_{i\in \mathbb{N}}$. 
That is,  $\U^i = \LC U^i_j \RC_{j\in J_i}$ covers $X$ and it is a refinement of $\U^{i+1}$ 
in that  there is a function $p_i:\U^i \goto \U^{i+1}$ so that $U \subset p_i(U)$.
The refinement function $p_i$ readily extends to a simplicial map $p_i : \mathcal{N(U}^i) \longrightarrow \mathcal{N(U}^{i+1})$, 
and  it is the basis (as would  any other such  map)  for our cascading algorithms.

\begin{rmk} \label{rmk:topo_compare}
Given a simplicial complex $K$, the coboundary operator $\delta:C^0( K ; \mathbb{R}) \goto C^1(K ; \mathbb{R})$ 
sends scalar functions on the vertices of $K$ to (antisymmetric) scalar functions on its edges.
Such coboundary maps set the stage for cohomology and Hodge theory.
Moreover, a simplicial map $p:K_1 \goto K_2$ induces linear maps $p^*: C^n(K_2;\mathbb{R} ) \goto C^n(K_1; \mathbb{R})$.
The fact that $p^*$ and $\delta$ commute is necessary for tracking the generators of the cohomology groups in persistent cohomology.
In similar fashion, taking the adjoint $\delta^\dag:C^1(K; \mathbb{R}) \goto C^0(K ;\mathbb{R})$ of  $\delta$, we can define the Hodge Laplacian operator 
$\triangle^0 = \delta^\dag\circ\delta:C^0(K;\mathbb{R})\goto C^0(K;\mathbb{R})$.
Unfortunately, $\delta^\dag$ and hence $\triangle^0$ do not (in general) commute with transformations induced by simplicial maps,
and therefore Laplacian eigenvalues/eigenvectors cannot be made to correspond under $p^*$, as is done in \cite{edelsbrunner2015persistent}. To see this, 
let  $K_1 = \mathcal{P}(\LC a,b,c \RC)\smallsetminus \{\emptyset\}$ (the full 2-simplex) $K_2 = \mathcal{P}(\LC A, B \RC)\smallsetminus\{\emptyset\}$, 
and let  $p:K_1 \goto K_2$ be the simplicial map  $p(a) = A$, $p(b) =  p(c) = B$. 
Let $f$ be the cocycle $f(A,B) = 1$. Then, $[\delta^\dag p^*f](a) = [p^*f](b,a) + [p^*f](c,a) = f(B,A) + f(B,A) = -2$.
On the other hand, $[p^*\delta^\dag f](a) = [\delta^\dag f](A) = f(B,A) = -1$, 
and thus the adjoint  
$\delta^\dag$ does not commute with $p^*$.
\end{rmk}

We next describe two procedures for building cover towers along with the edge weights and partitions of unity respectively used to define the pertinent Laplacian operators and approximate eigenfunctions.
The subsequent sections  describe the cascading algorithms aimed at addressing the problem described in Remark \ref{rmk:topo_compare}.

\begin{example}[Cover Tree Graphs] \label{ex:graph_CT} \rm
    The first construction builds a cover tower via a cover tree   
    \cite{Cover_Trees_1, Cover_Trees_2}, which is a data structure originally designed for efficient neighborhood query algorithms.
    \begin{defn} \label{def:CT} \rm
    Consider a dataset $X$ within a metric space $(\X,d)$.
    A cover tree $T$ on $X$ is a directed tree with nodes labeled by points in $X$ and edges indicating assigned nearby children.
    All nodes at fixed depth $i$ determine a level set $C_i$ of the tree.
    A dyadic cover tree \cite{Cover_Trees_1} has the following properties:

    \begin{itemize}
        \item (nesting) $C_i \subset C_{i-1}$.
        \item (cover) for each $x \in C_i$ and each $j > i$ there is a unique $y \in C_j$ so that $y$ is the grandparent of $x$, satisfying $d(x,y) < 2^{j+1} - 2^{i} = \sum_{k=i}^j 2^k$. Denote this relationship as $y = GP_j(x)$.
        \item (separation) For distinct $x, x' \in C_j$ in the same level, $d(x,x') > 2^j$.
    \end{itemize}
    \end{defn}

    The cover and nesting properties guarantee that $\U^i = \LC B(x,R\times2^{i+1}) : x \in C_i \RC$ defines a cover tower of $X$ for any fixed ratio $R \geq 1$, and the unique grandparent assignment automatically yields a refinement  function.
    The separation property controls the growth rate of $\Ln C_i \Rn$ and yields evenly spaced covers.

    To define a partition of unity, consider a radial kernel $k_h(x,y) = K(d(x,y)/h)$ with $\textrm{supp}(K) = [0,1]$ (e.g., $K(r) = \max(1-r,0)$).
     Thus $\textrm{supp}(k_h(x,\cdot)) = \overline{B(x, r)}$ and for \hbox{$r_i=R\times2^i$} from the cover tree construction,
    then $\phi^i_x(y) = \frac{k_{r_i}(x,y) }{ \sum_{x \in C_i} k_{r_i}(x,y)}$
    defines a partition of unity dominated by $\U^i$ for each scale $i$.
    In similar fashion,   $W^i_{kj} = \Ln GP_j^{-1}(x_k) \Rn \Ln GP_j^{-1}(x_j) \Rn k_{2*r_i}(x_k,x_j)$ defines density-aware edge weights,
    wherein the domain of $ GP_j $ is $ \cup_{i < j} C_i \subset  X$.
\end{example}

\begin{example}[TDA Mapper Graphs] \label{ex:graph_mapper} \rm
    The mapper method \cite{mapper} generates graphs from a dataset $X$ and a filter function $f:X \goto \R$ by splitting the filter-space $[min(f(X)),max(f(X))]$ into a number of overlapping intervals $I_k$, and dividing each $f^{-1}(I_k)$ into clusters $X_i$.
    The collection of clusters $X_i$ is a cover of $X$ and the 1-skeleton of its nerve is the mapper graph.
    Refining the interval cover $\{I_k\}$ also refines its pullback under $f$, but not necessarily the clusters $X_i$; without refinement from one level to the next, we do not necessarily obtain a cover tower. 
    Moreover, we lack parent maps for cascade; instead we average over each cover set the approximate eigenfunctions from the coarser cover.
    This approach is reasonable to compute for mapper graphs.

    To yield a partition of unity, a Gaussian (or other) density $g_i$ is fitted to each partial cluster $X_i$ and the collection of densities is renormalized to obtain $\phi_i(x) = g_i(x)/\sum_j g_j(x)$.
    Edge weights are defined according to intersection size $W_{ij} = \Ln X_i \cap X_j \Rn$.
    The filter function, number of intervals, degree of overlap, and cluster method are all parameters of the mapper construction.
    The cluster method must be flexible with the number of clusters, such as thresholded hierarchical clustering (here we use the recommended single linkage).
\end{example}

\subsection{Graph Laplacians}
A given similarity matrix $W_{kj}$ defines multiple graph Laplacians.
\begin{defn} \label{def:laplacians}
    Let $G = (V,E)$ be a graph with edge-weight matrix $W$ and
    $D$ the diagonal degree matrix with $D_{ii} = \sum_j W_{ij}$.
    The unnormalized graph Laplacian $L$, normalized symmetric graph Laplacian $L_{sym}$, and random walk graph Laplacian $L_{rw}$ are defined  as
    \begin{align}
        L &= D - W, \\
        L_{sym} &=  I - D^{-1/2}WD^{-1/2}, \\
        L_{rw} &=    I - D^{-1}W.
    \end{align}
\end{defn}

The unnormalized graph Laplacian $L$ has poor spectral convergence (as  $\#(V)\goto \infty$) properties   \cite{Graph_Lap_Converge}. 
Thus, we focus on the normalized symmetric and random walk graph Laplacians together, as they are similar matrices ($L_{rw}v = \lambda v$ if and only if $L_{sym}D^{1/2}v = \lambda D^{1/2}v$).
To contrast, $L_{sym}$ is symmetric and therefore convenient for eigenvector determination (e.g., via Lanczos \cite{lanczos} or LOBPCG \cite{lobpcg} methods), while $L_{rw}$ always has the constant vector as a null eigenvector.
Finally, $L_{rw}$ is a Hodge Laplacian $\delta_w^\dag \circ \delta_w$ as described in \cite{Combo_Laplace},  for a compatible weighted coboundary operator.

\subsection{Eigenvector Cascading}
Here we present our two algorithms, dubbed first and second eigenvector cascade. 
In general, these methods apply to a series of matrices $L_i$ associated to graphs $G_i$ with simplicial maps $p_i$ between them.
Since the mappings $p_i$ should be geometrically defined, the matrices $L_i$ should reflect geometry underlying the graphs, such as a Laplacian operator;
thus, we stress that the cascade methods are designed to align geometric features across scales.

First cascade, shown in Alg. \ref{alg:cascade1}, is designed to accelerate eigenvector determination by taking advantage of the economical initial guesses provided by $p_i^*$.
Only as many eigenvectors $v_j^i$ as are initially determined ($m$) can be cascaded, though additional guess vectors $u_j^i$ can be added, e.g. chosen randomly followed by Gram-Schmidt.

\begin{algorithm}[ht] \label{alg:cascade1}
\SetAlgoLined
\DontPrintSemicolon
\KwData{Simplicial maps $p_i$,
Matrix operators $L_i$, \\
Eigensolver $F(L,V)$ for $L$ and initial guess V, \\
A starting  state $\LC u_j^N \RC_{j=1}^m$ at the coarsest scale.}
\KwResult{Eigenvectors $\LC v_j^i \RC_j$, eigenvalues $\LC \lambda_j^i \RC_j$ and initial states $\LC u_j^i \RC_j$ for each scale $i$.}
\Begin{
 	$\LP \LC v_j^N \RC_{j=1}^{m}, \LC \lambda_j^N \RC_{j=1}^{m} \RP = F\LP L_N, \LC u_j^N \RC_{j=1}^{m} \RP$ \;
 	\For{$i$ from $N-1$ to $n$}{
 	    $\LC u_j^i \RC_{j=1}^{m} = p_i^*\LP\LC v_j^{i+1} \RC_{j=1}^{m}\RP$ \;
 	    $\LP \LC v_j^i \RC_{j=1}^{m}, \LC \lambda_j^i \RC_{j=1}^{m} \RP = F\LP L_i,\LC u_j^i \RC_{j=1}^{m}\RP$
}   }
\caption{First Cascade}
\end{algorithm}

Second cascade, shown in Alg. \ref{alg:cascade2}, forwards the vertex relationships defined by $p_i$ to the eigenvectors via repeated projection onto the eigenspaces defined by first cascade.
These eigenspaces are organized by partitions $\LC P_i[k] \RC_{k=1}^{N_i}$ of the eigenvector indices for each scale $i$. 
The need for second cascade is exemplified when eigenvalues are clustered near each other; in this case, individual eigenvectors are unstable under small perturbation (c.f. \cite{EV_perturb} or \cite{ES_proj_perturb} theorem 2.1).
Thus, second cascade focuses on the stable portion---their span---and uses projection to create relatable bases, i.e. minimizing the errors $\Ln p_i^*(w^{i+1}_j) - w^i_j \Rn$).
Indeed, though $w^i_j$ need not be eigenvectors, the spans $\LC w^i_j \RC_{j \in P_i[k]}$ and $\LC v^i_j \RC_{j \in P_i[k]}$ for each eigenvalue cluster $P_i[k]$ are identical.

In practice, the output of Alg. \ref{alg:cascade1} is sufficient for Alg. \ref{alg:cascade2} and are performed in tandem as \emph{double cascade} for a collection of Laplacian operators $\LC L_i \RC$.
Toward adaptive eigenvalue clustering, we group eigenvalues if they are within some fixed multiple of machine error or if their ratio is smaller than a fixed value.

\begin{algorithm}[!ht] \label{alg:cascade2}
\DontPrintSemicolon
\KwData{Paired eigenvectors $\LC v_j^i \RC_{j=1}^m$ and eigenvalues $\LC \lambda_j^i \RC_{j=1}^m$ (in increasing order) for each scale $n \leq i \leq N$, parent maps $p_i$, and small grouping thresholds $\delta$ and $\epsilon$.}
\KwResult{A collection of new basis vectors $\LC w_j^i \RC_j$ for each scale $i$, and proposed eigenspace groupings, $P_i$ which partition $\LC 1,...,m \RC$.}
\Begin{
    $\LC w_j^N \RC = \LC v_j^N \RC$

    \For{$i$ from $N-1$ to $n$}{
        $P_i = \LP \LP 1 \RP \RP$ \;
        $N_i = 1$ \;
        \For{$j$ from $2$ to $m$}{
            \If{$\Ln \lambda_{j-1}^i - \lambda_{j}^i \Rn < \delta$ or $\lambda_{j}^i/\lambda_{j-1}^i < 1+\epsilon$}{
                Append $j$ to $P_i[N_i]$
            } \Else {
                Append $\LP j \RP$ to $P_i$ \;
                $N_i = N_i + 1$
        }   }
        \For{$k$ from $1$ to $N_i$:}{
            $V = \textrm{span}\LP \LC v_j^i\RC_{i \in P_i[k]} \RP$ \;
            \For{$j$ in $P_i[k]$}{
                $w_j^i = Proj_V\LP p_i^*(w_j^{i+1})\RP$
}   }   }   }
\caption{Second Cascade}
\end{algorithm}

\section{Experimental Results} \label{S:examples}

\subsection{Speeding Up Laplacian Eigenvector Computations with First Cascade} \label{SS:time}

\begin{example} \label{ex:timers} \rm
    Here we examine the computational benefit of first cascade (Alg. \ref{alg:cascade1}).
    We do not expect to establish state-of-the-art results here, but rather hope to re-affirm results such as \cite{F1_cascade} and \cite{ACE_Cascade} in a different context where the graphs and Laplacian operators are collapsed in a different fashion.
    We compare the run-time of LOBPCG via paired trials on several generated datasets both with and without cascading initial states.
    The trial results are given in Table \ref{tab:timing}.
    Cascading always yields a noticeable speedup at each scale and in most cases total cascade time undercuts na\"{i}ve LOBPCG performed only at the finest cover level (the most complex eigensystem problem).

    \begin{table}[htb!]
        \centering
        \begin{tabular}{c|ccc}
    dataset	& single cascade &	lobpcg final & lobpcg full \\ \hline
    PIN\_1  & \textbf{28.84}   & 55.36  & 86.55    \\
    PIN\_2  & \textbf{30.87}   & 67.49  & 74.5     \\
    CANTOR  & \textbf{3.53}    & 4.59   & 12.9     \\
    CARPET  & \textbf{117.06}  & 125.34 & 152.66   \\
    SPHERE  & 40.02	  & \textbf{37.19}  & 44.48    \\
    BOXTREE & 60.89	  & \textbf{56.08}  & 68.03    \\
    EDGES   & \textbf{305.09}  & 346.58 & 419.07   \\
        \end{tabular}
        \caption{Time in seconds to find 100 (smallest) Laplacian eigenvectors for generated datasets at multiple scales.
        The cover tree construction of Ex. \ref{ex:graph_CT} is applied to each dataset.
        Shortest times are shown in bold.
        Each dataset has roughly 15,000 points and intrinsic dimension 2, except EDGES contains 100,000 points and CANTOR and CARPET sample fractals with respective Hausdorff dimensions $1 < \log_3(4),\log_3(8) < 2$.}
        \label{tab:timing}
    \end{table}


    The CANTOR dataset attains the largest speedup.
    As a sample from the Cantor square, many connected components associate to $\mathsf{ker}(L_{rw})$;
    thus, convergence is attained immediately once 100 connected components are fully resolved (ie, the initial guesses are already eigenvectors).
    The resulting speedup is recorded in Table \ref{tab:cantor}.

    \begin{table}[htb!]
        \centering
        \begin{tabular}{c|c c}
        depth & non-cascade & single cascade \\ \hline
        8	  & 0.42             & 0.17 \\
        9	  & 0.79             & 0.20 \\
        10	  & 4.44             & 2.33 \\
        11	  & 2.66             & 0.60 \\
        12	  & 4.59             & 0.24 \\
        \end{tabular}
        \caption{Time in seconds to perform LOBPCG with and without cascading initial states for the CANTOR dataset at various depths in the cover tree. There is a large drop in cascade time after depth 10, for which $\textrm{dim(Null)}(L_{rw}) \geq 100$.}
        \label{tab:cantor}
    \end{table}
\end{example}

\subsection{Double Cascade of Laplacian Eigenspaces} \label{SS:double_cascade}

Here we demonstrate using double cascade (Algs. \ref{alg:cascade1} and \ref{alg:cascade2} in tandem) to obtain consistent eigenspace basis representation for symmetric shapes, which have repeated (in the limit) or clustered (in data) eigenvalues.
Since these are low-dimensional examples, we can view the approximate eigenfunctions defined by the nerve map (see Fig. \ref{fig:NL_vignette}) as a common space to compare across scales.
To enforce the effects of symmetry, the example datasets are given as grids of points, though the cover trees do not always reflect this.
Though symmetry elucidates the effects of double cascade, it is helpful for any clustered eigenvalues; even for isolated eigenvalues, double cascade enforces consistency in eigenvector sign.

\begin{example} \label{ex:PIN_1} \rm
    The PIN\_1 dataset consists of points evenly dispersed within three triangles arranged to be invariant under 120 degree rotation.
    Since PIN\_1 is connected, the first (null) eigenvector is ignored.
    Heat maps of the eigenfunctions derived from single and double cascade are shown in Fig. \ref{fig:pinwheel}.
    Comparing these cases reveals an additional benefit to second cascade:
    In addition to consistency, the basis vectors originate at and reflect a coarse scale;

    \begin{figure}[ht!]
        \begin{center}
            Single Cascade \\
            \includegraphics[scale = 0.25]{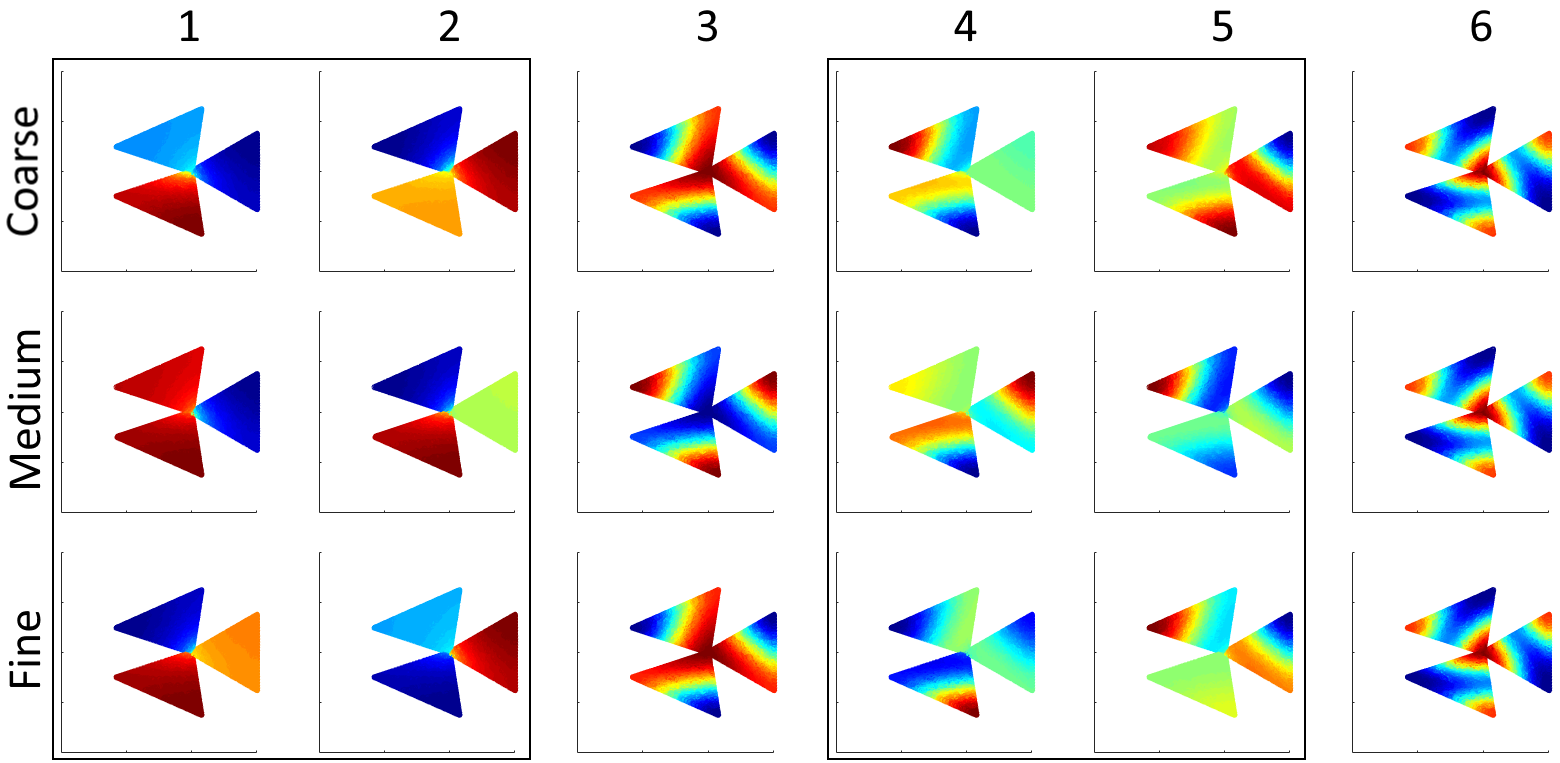} \\
            Double Cascade \\
            \includegraphics[scale = 0.25]{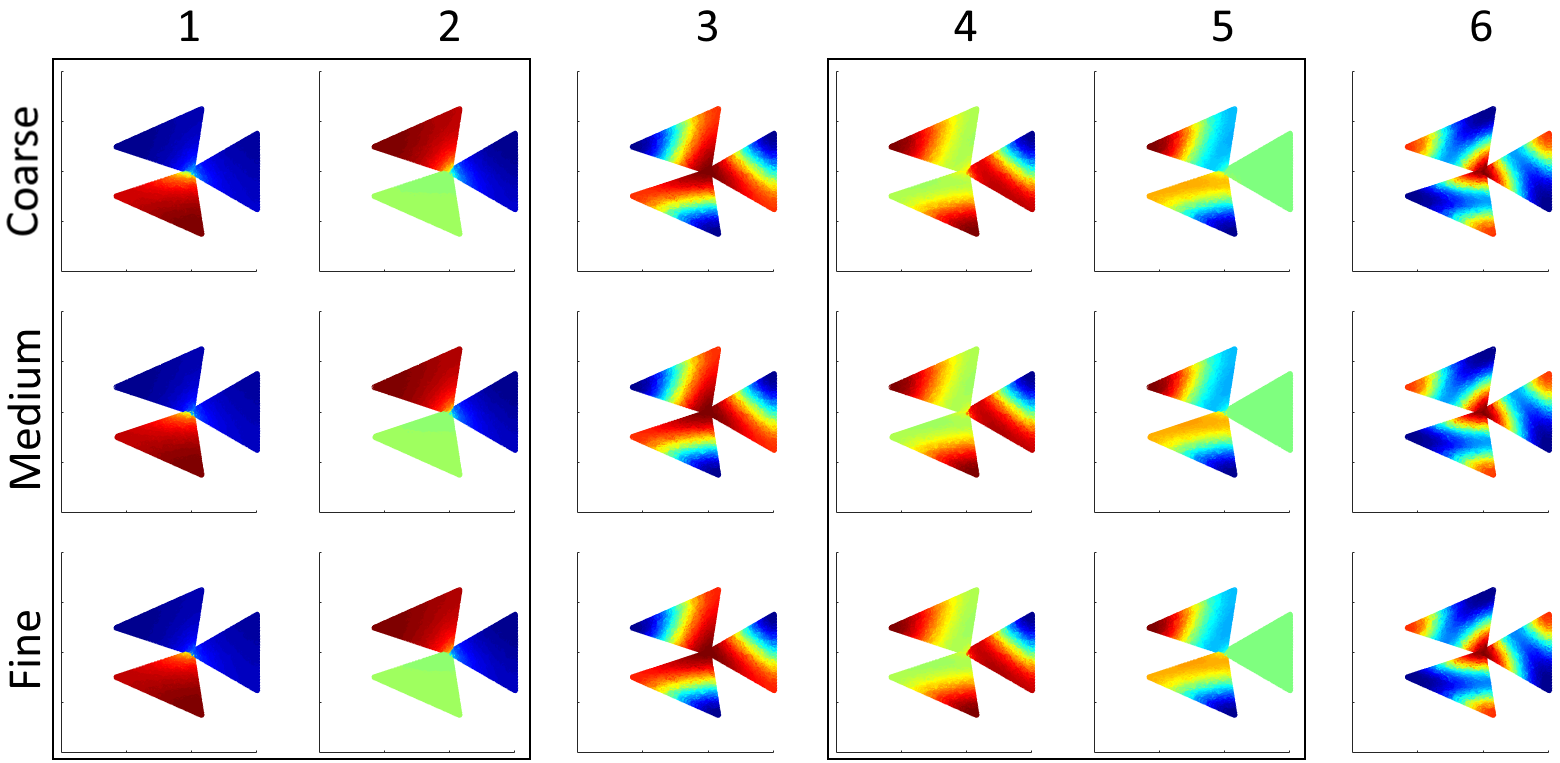}
        \end{center}
        \caption{Nontrivial Laplacian eigenfunction heat maps for the PIN\_1 dataset at multiple scales for both single and double cascade.
        Each row is a choice of scale while each column is an eigenfunction index.
        Each eigenspace has a large bounding box.}
        \label{fig:pinwheel}
    \end{figure}
\end{example}

\begin{example} \label{ex:SPHERE} \rm
    The SPHERE dataset consists of six tangent geodesic annuli on a 2-sphere, and exhibits cubical symmetries.
    Since the dataspace is connected at any positive resolution, the first (null) eigenvector is ignored.
    Heat maps of the eigenfunctions derived from single and double cascade are shown in Fig. \ref{fig:sphere};
    In this case, eigenspace convergence (and symmetries) is evident with double cascade and obscure without it.

    \begin{figure}[htb!]
        \begin{center}
            Single Cascade \\
            \includegraphics[scale = 0.25]{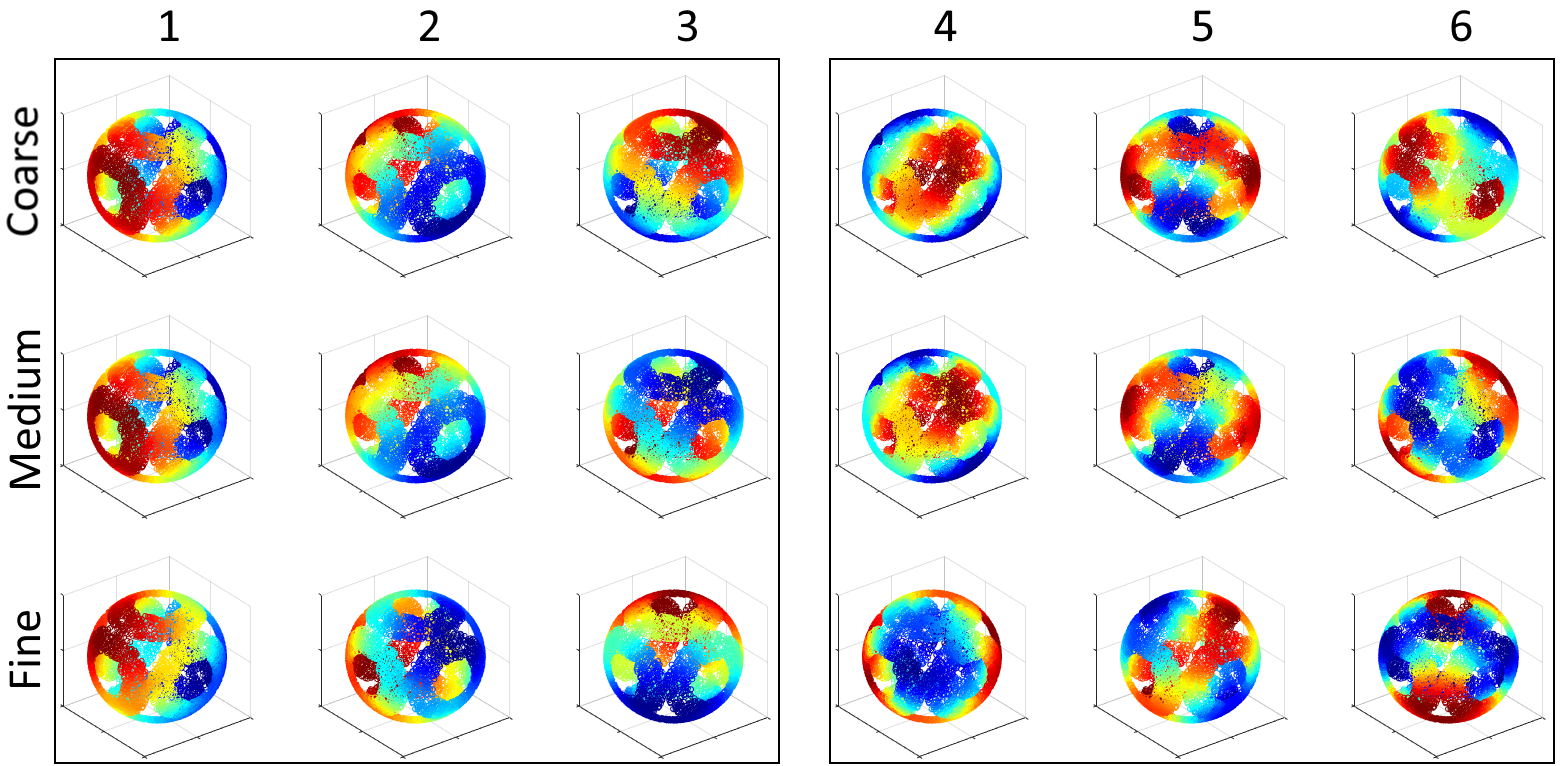} \\
            Double Cascade \\
            \includegraphics[scale = 0.25]{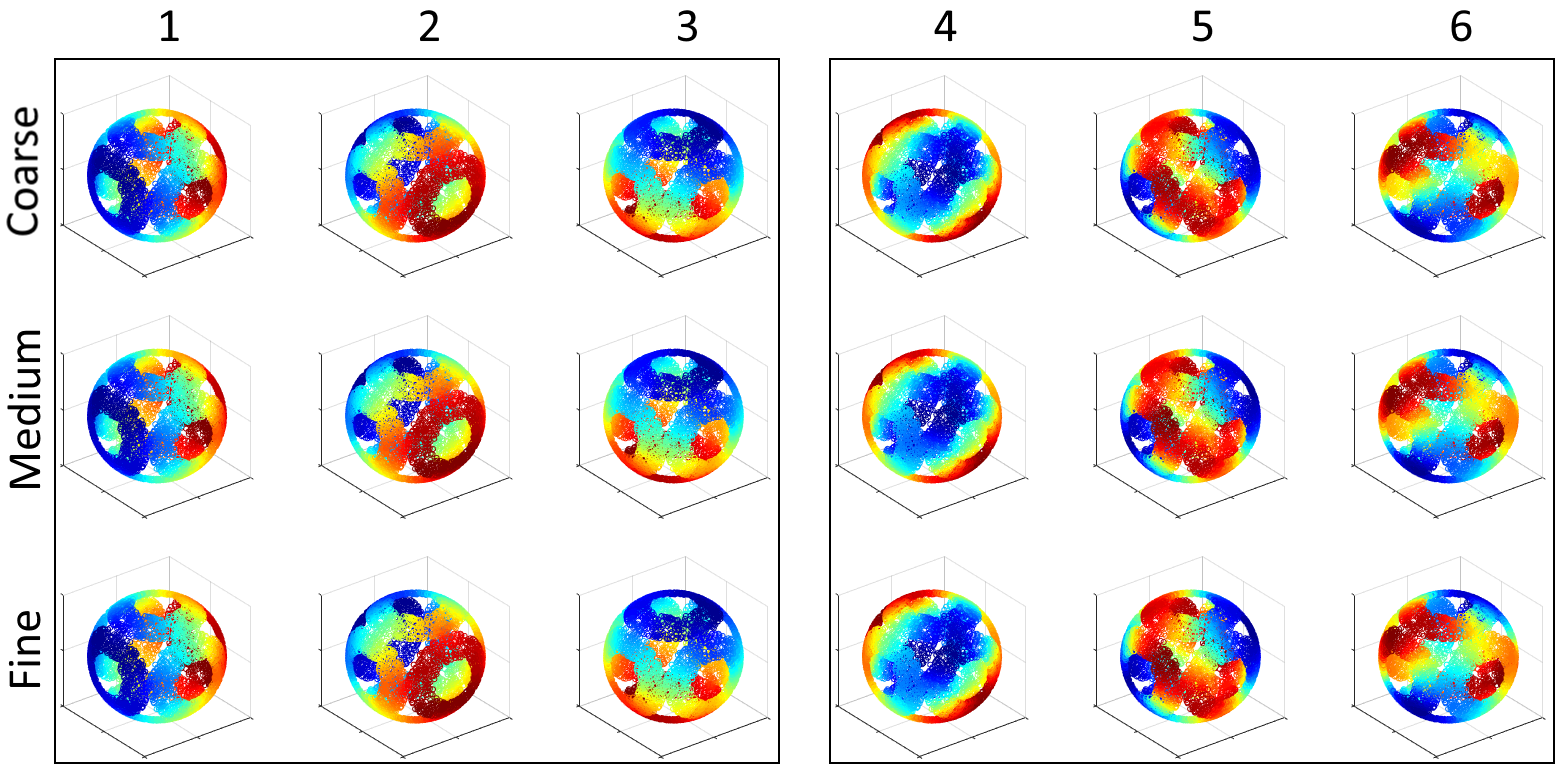}
        \end{center}
        \caption{Nontrivial Laplacian eigenfunction heat maps for the SPHERE dataset at multiple scales for both single and double cascade.
        Each row is a choice of scale while each column is an eigenfunction index.
        Each eigenspace has a large bounding box.}
        \label{fig:sphere}
    \end{figure}
\end{example}

\begin{example} \label{ex:CANTOR} \rm
    The CANTOR dataset consists of the seventh iteration of the Cantor square fractal.
    With a finite approximation, the Cantor square emulates a dataset with many connected components which coalesce at different scales.
    Heat maps of the eigenfunctions derived from single and double cascade are shown in Fig. \ref{fig:cantor}.
    Despite having equal (null) eigenvalues, double cascade naturally orders the eigenvector bases by separation scale;

    \begin{figure}[htb!]
        \begin{center}
            Single Cascade \\
            \includegraphics[scale = 0.25]{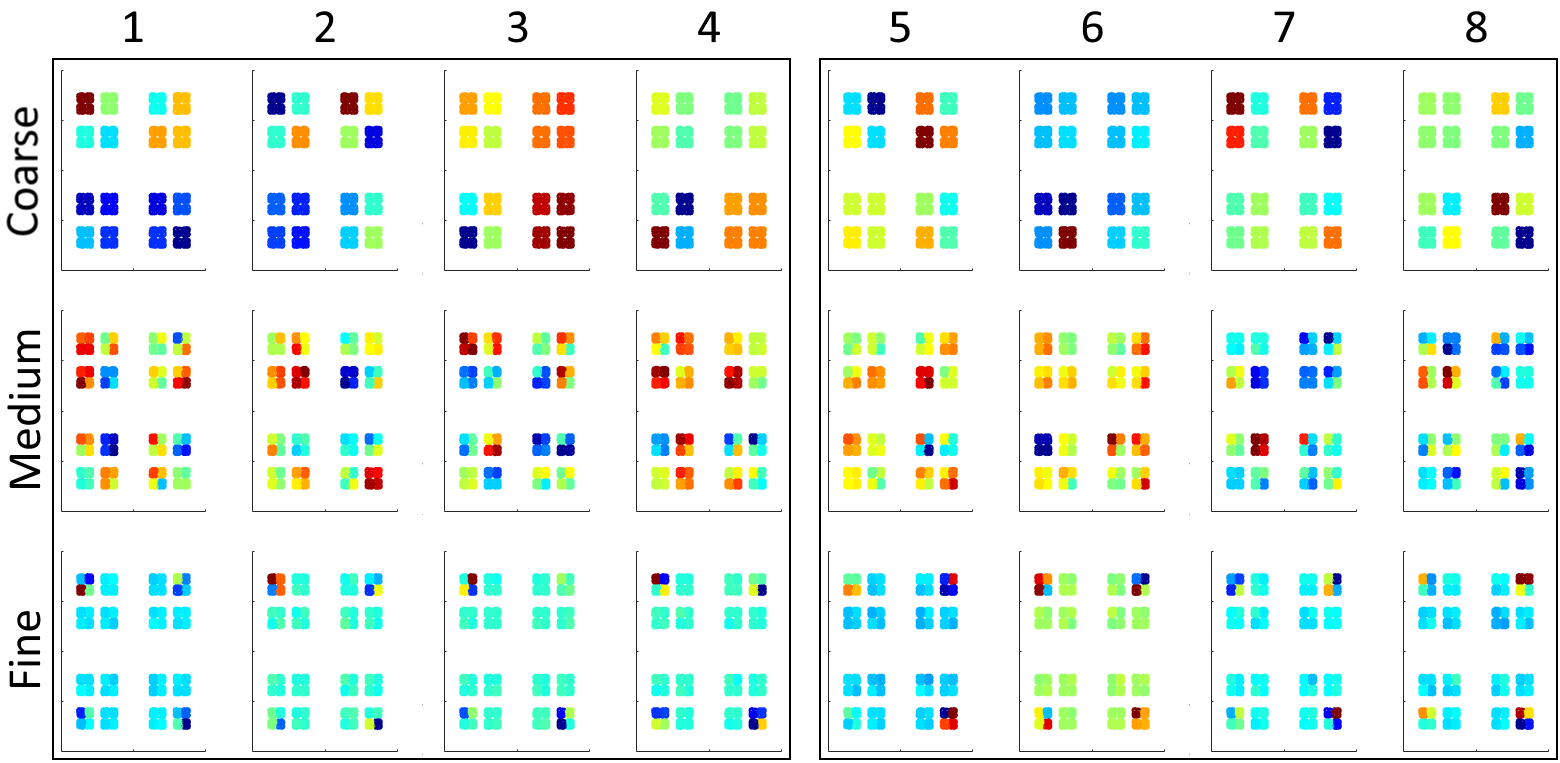} \\
            Double Cascade \\
            \includegraphics[scale = 0.25]{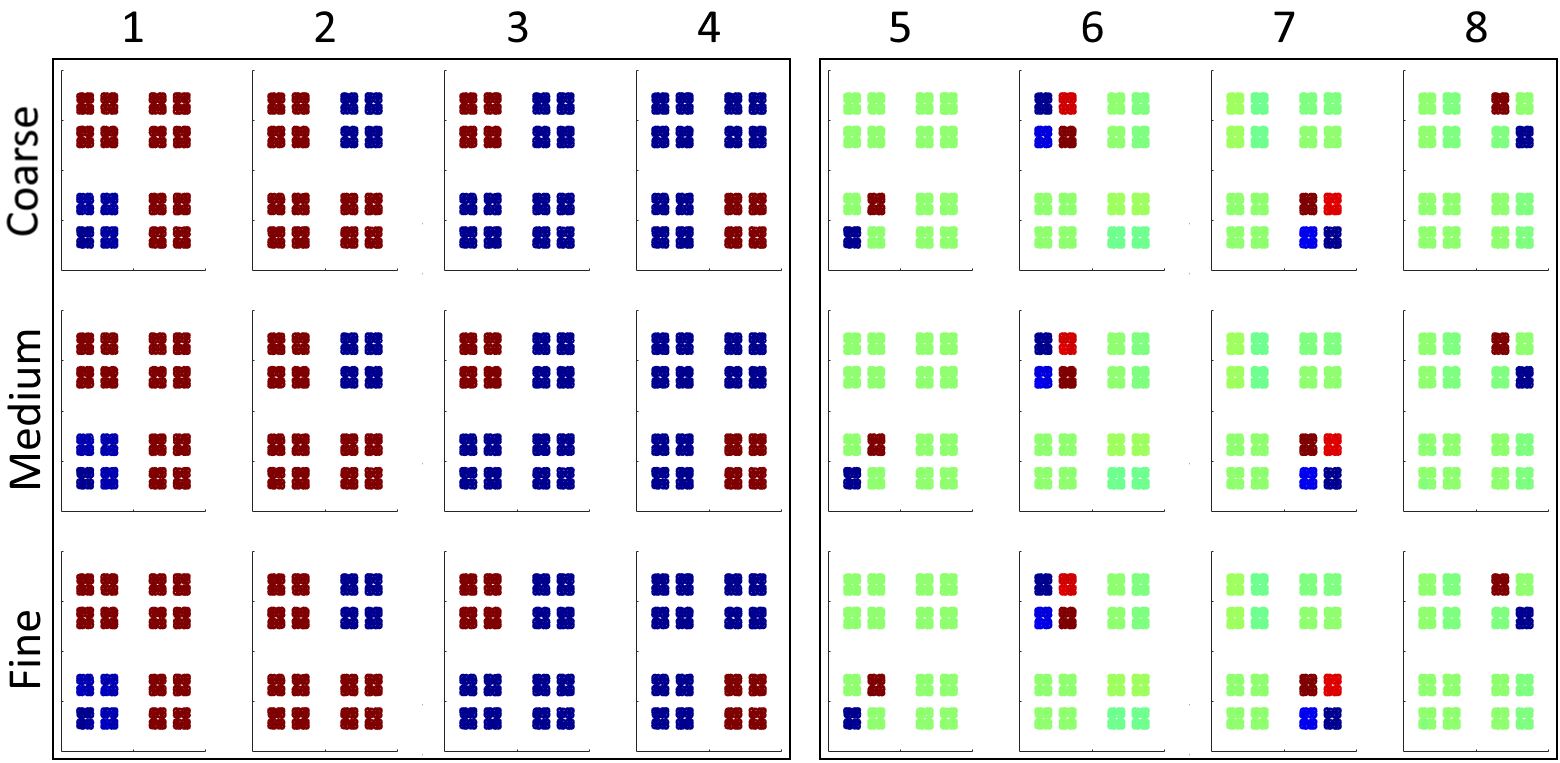}
        \end{center}
        \caption{Laplacian eigenfunction heat maps for the CANTOR dataset at multiple scales for both single and double cascade.
        Each row is a choice of scale while each column is an eigenfunction index.
        Each eigenspace has a large bounding box.
        Due to the fractal nature of the dataset, the first four eigenvectors are grouped, the fifth through $16^{th}$ are grouped, and so on in powers of four.}
        \label{fig:cantor}
    \end{figure}
\end{example}

\subsection{Spectral Persistence in Multiscale Mapper} \label{SS:mapper}

Next we investigate the progression of Laplacian eigenfunctions for mapper graphs built at ramping scale via double cascade for each of three real-world datasets.
In this case, we are most interested in establishing basis correspondence across resolution to analyze the dataset geometry across scales, specifically via weakly connected components or flares highlighted by Laplacian eigenvectors with small eigenvalue.
Due to dataset symmetries, such regions can be very difficult to track in the mapper graphs;
double cascade accomplishes this task when graph embeddings fail.
Note that even na\"ive spectral graph embeddings may fail in this task since symmetry causes instability in individual eigenvectors (see subsection \ref{SS:double_cascade}).
Similar to \v{C}ech persistent homology \cite{TDA}, we expect that important eigenspaces will be represented at a broad range of scales with spurious features relegated to extreme or specific scales.

To avoid spurious connected components in the graph and thereby simplify presentation, mapper clusters with only one point are removed.
Such singletons generally represent noise or the cut-off tip of a flare; none were found to persist.
The following datasets have high intrinsic dimension and nonlinear shape, making the approximate eigenfunctions impossible to view.
Instead, the eigenvectors defined on the mapper graphs' vertices will be compared directly via Fruchterman Reingold force-directed layout \cite{FR_graph_force} in 2D.
Consequently, double cascade is preferred to properly correspond the eigenvector features from one scale to the next, even over direct inspection.

\begin{example} \label{ex:diabetes} \rm
    Consider the diabetes dataset \cite{Diabetes_Data} of diagnostic measurements from patients with chemical or overt diabetes.
    This dataset was investigated in the initial mapper article \cite{mapper}.
    As in the original analysis, eccentricity is chosen as the filter function.
    The resulting mapper graphs with eigenvector heat maps are plotted in Fig. \ref{fig:diabetes}.
    Sequential eigenfunctions show a great deal of similarity despite the changing underlying graphs;
    moreover, each eigenvector concentrates within a flare in the mapper graph indicating a region of potential importance.

    \begin{figure}[htb!]
        \begin{center}
            \includegraphics[scale=0.3]{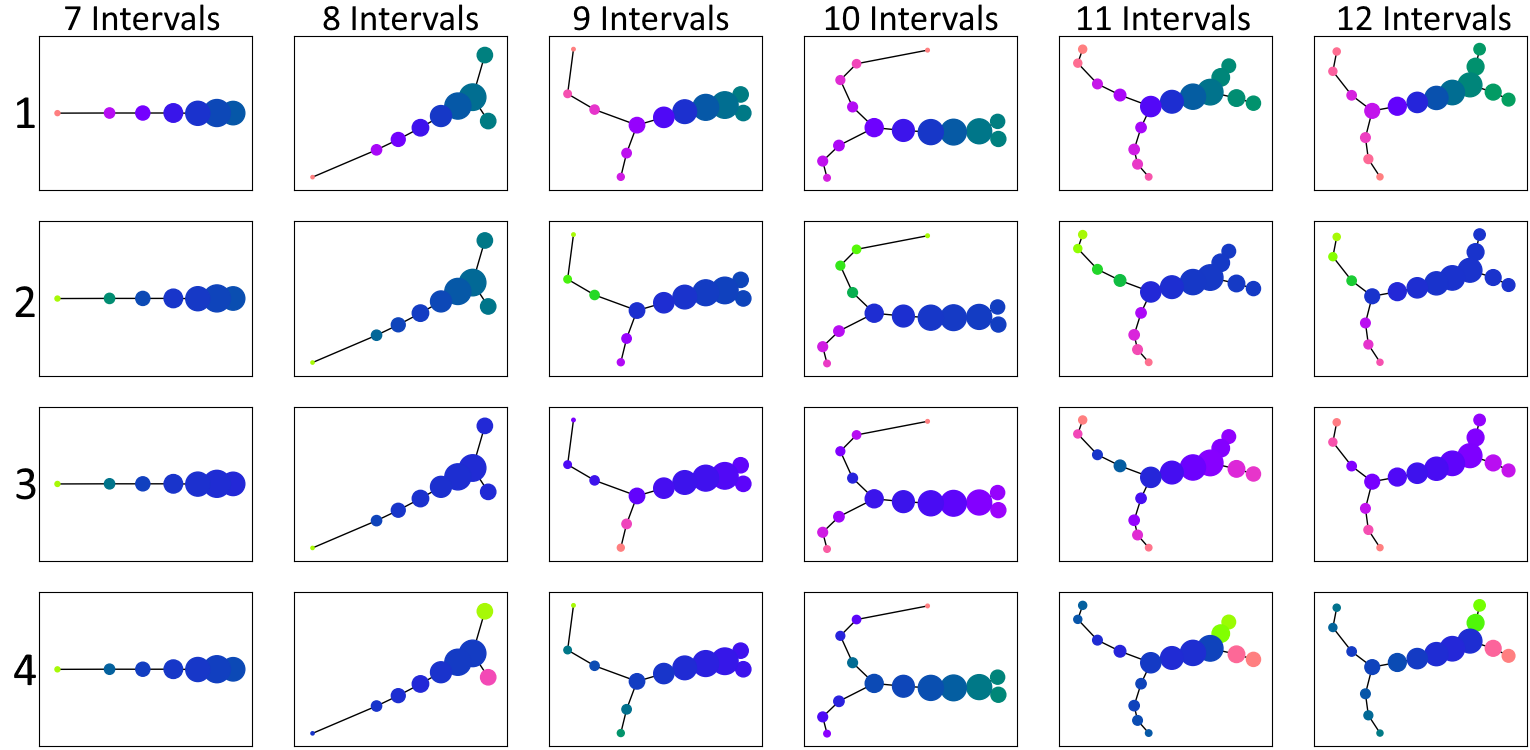}
        \end{center}
        \caption{Eigenvectors for the mapper graphs obtained from the diabetes dataset \cite{Diabetes_Data}.
        Each row indicates an eigenvector index, while each column represents a particular scale (number of intervals).}
        \label{fig:diabetes}
    \end{figure}
\end{example}

\begin{example} \label{ex:iris} \rm
    Next consider Fisher's iris dataset \cite{Fisher_Iris, Andersen_Iris}.
    The filter function is a Gaussian density estimate, wherein the bandwidth is the average distance to the $10^{th}$ nearest neighbor.
    The resulting mapper graphs with eigenvector heat maps are plotted in Fig. \ref{fig:iris}.
    Here Eigenvector correspondence is a great help in relating the mapper graphs, though the finest resolution (30 intervals) breaks into separate components.

    \begin{figure}[htb!]
        \begin{center}
            \includegraphics[scale=0.3]{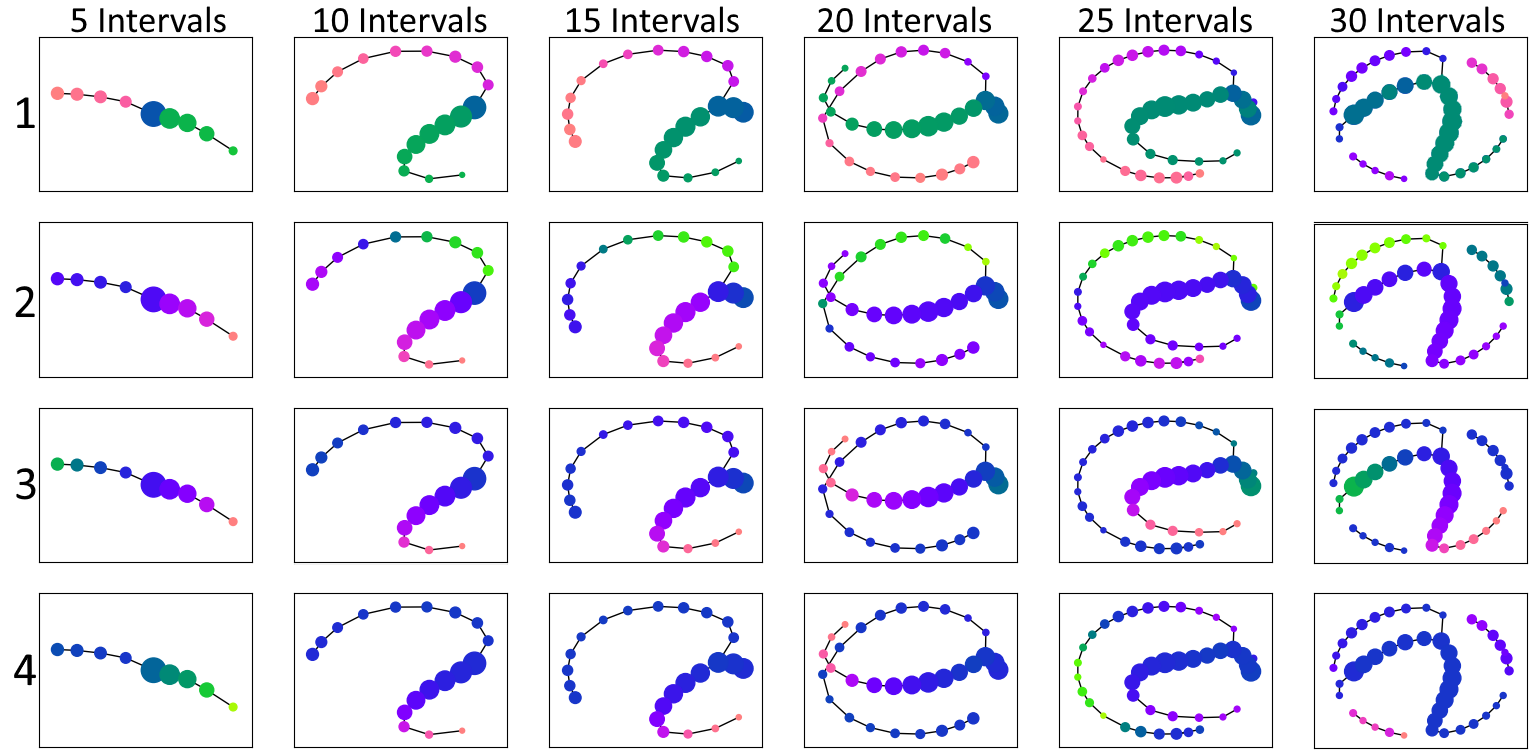}
        \end{center}
        \caption{Eigenvectors for the mapper graphs obtained from the iris dataset \cite{Fisher_Iris,Andersen_Iris}.
        Each row indicates an eigenvector index, while each column represents a particular scale (number of intervals).}
        \label{fig:iris}
    \end{figure}
\end{example}

\begin{example} \label{ex:images} \rm
    Consider the dataset of high-contrast 3x3 grayscale natural image patches, sampled randomly from the Van Hateren dataset \cite{Van_Hateren}.
    Similar data is investigated in \cite{TDA_Natural_Images} with a much larger sample size.
    For mapper, the filter function is a Gaussian density estimate with bandwidth chosen as the mean distance to the $10^{th}$ nearest neighbor.
    The resulting mapper graphs and eigenvectors are plotted in Fig. \ref{fig:VH}.
    The four flares represent connected components of the low density region, and support the 3-circle model presented in \cite{TDA_Natural_Images}.
    Though the graph layout permutes these 4 flares, the eigenfunctions relate very stably from 8 intervals to 14 intervals, evidenced by comparing the associated eigenfunctions.

    \begin{figure}[htb!]
        \begin{center}
            \includegraphics[scale=0.3]{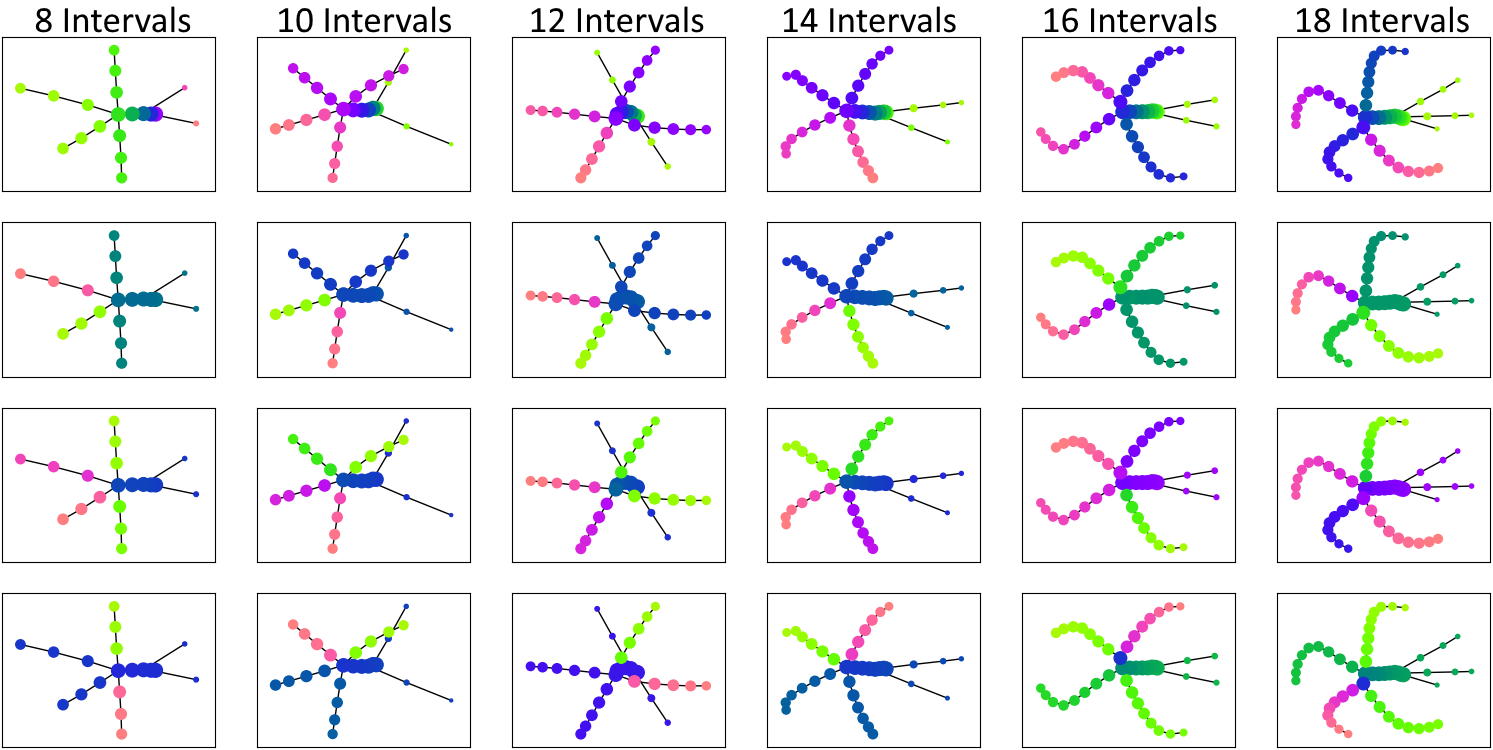}
        \end{center}
        \caption{Eigenvectors for the mapper graphs obtained from high-contrast 3x3 natural image patches \cite{TDA_Natural_Images,Van_Hateren}.
        Each row indicates an eigenvector index, while each column represents a particular scale (number of intervals).}
        \label{fig:VH}
    \end{figure}
\end{example}

\section{Discussion and Conclusion} \label{S:CD}
    In this work we view the progression of Laplacian eigenvectors obtained from a sequence of graphs consistently across scale.
    A cover tower framework is proposed to construct graphs from a dataset at multiple scales, instantiated via cover trees and mapper.
    We discuss why existing methods (the persistent homology of a self map \cite{edelsbrunner2015persistent}) do not apply here, and thus propose two cascade methods to aid in corresponding eigenspace bases.
    Our experiments show that first cascade accelerates LOBPCG calculation of the eigenvectors and that second cascade elucidates   relationships between eigenvectors at different scales.

    Next, we used the consistent bases from double cascade to properly relate mapper graphs built at multiple scales.
    In this context, eigenvector cascade can be used to make spectral visualization (dimensionality reduction) consistent across scale and aids in using mapper for dataset exploration and the identification of unusual regions.
    This story parallels that of persistent homology;
    specifically, important dataset features persist for a broad range of scales, while others are transient.
    That said, we do not track loops or voids (as in co-homology), but weakly-connected components identified by the support of Laplacian eigenvectors with small eigenvalue.

    Companion work develops a framework for weighted graphs and graph collapse underlying the experiments presented here.
    There, the notion of graph collapse is linked to simplicial maps;
    This analysis invites a more general notion of (weighted) simplicial collapse as a generalization to the parent maps in cascading.
    Moreover, such simplicial maps may also be used to cascade eigenvectors of weighted $n$-Laplacians (c.f. \cite{Combo_Laplace}) and thereby track a richer family of geometric features.

    Our investigation leads us to consider a notion of persistent Laplacian eigenspaces (or PLES) as a geometric extension to persistent cohomology.
    In particular, PLES can be used to track cohomological features as null eigenvectors.
    Since Laplacian operators do not commute with simplicial maps, as in traditional persistence, 
    this leads one to the problem of geometrically meaningful tracking of the eigenfunctions themselves.
    In this context, eigenvector cascade provides tools to handle the rise and fall of eigenvalues, in particular by consistently tracking eigenvectors through regions of clustered eigenvalues.
    We see this approach as useful for more than multiscale analysis.
    So long as partial correspondence can be found between graph vertices, double cascade can be applied to track persistent geometric features in time-varying graphs including language graphs, dynamic sensor networks, social networks, and so on.



\bibliography{NL}{}
\bibliographystyle{plain}
\end{document}